\documentstyle{article}

\newtheorem{theorem}{Theorem}

\newtheorem{definition}[theorem]{Definition}
\newtheorem{example}[theorem]{Example}
\newtheorem{corollary}[theorem]{Corollary}

\def\QED{\quad\blackslug\lower 8.5pt\null}

\begin{document}

\begin{center}
{\Large \bf ON SOME METHODS OF CONSTRUCTION} 

\vspace*{3mm}
{\Large \bf OF INVARIANT NORMALIZATIONS} 

\vspace*{3mm}
{\Large \bf OF LIGHTLIKE  HYPERSURFACES}\footnote{{\bf 
1991 MS classification}: 53B25, 53B20, 53B21, 53B30.

\hspace*{1mm} {\bf Keywords and phrases:} 
Pseudo-Riemannian manifold,  Lorentzian signature, lightlike 
hypersurface, invariant normalization, affine connection, 
  isotropic geodesics, singular point, isotropic 
sectional curvature.}

\vspace*{3mm}
{\large M.A. Akivis and  V.V. Goldberg}
\end{center}

\vspace*{5mm}

{\footnotesize{\it Abstract}. The authors study the geometry of 
lightlike hypersurfaces on pseudo-Riemannian manifolds $(M, g)$ 
of Lorentzian signature. 
Such hypersurfaces are of interest in general 
relativity since they can be models of different types of 
physical horizons. For a lightlike hypersurface 
$V \subset (M, g)$ of general type and for  some special 
lightlike hypersurfaces (namely, for 
totally umbilical and belonging to a manifold $(M, g)$ 
of constant curvature),  in a third-order neighborhood 
of a point $x \in V$, the authors construct  invariant 
normalizations  intrinsically connected  with 
the geometry of $V$ and investigate affine connections induced by 
these normalizations. For this construction, they used 
relative and absolute invariants defined by the first and second 
fundamental forms of $V$. The authors show that if 
$\dim M = 4$, their methods allow to construct three 
invariant normalizations and affine connections 
 intrinsically connected with  the geometry of $V$. 
Such a construction  is given 
in the present paper for the first time. 
The authors also consider the fibration of isotropic geodesics 
of $V$ and investigate their singular points and singular 
submanifolds. 
}

\vspace*{5mm}

\setcounter{equation}{0}

\vspace*{5mm}


\setcounter{equation}{0}

{\bf 0. Introduction.}  
The lightlike hypersurfaces $V$ of a  pseudo-Riemannian manifold 
$(M, g)$ of Lorentzian signature produce models of horizons of 
different types in general relativity. This is the reason 
they were studied intensively by geometers and physicists 
(see the books \cite{DB96}, \cite{Ku96}, 
\cite{HE73}, and \cite{HP95} as well as 
many papers quoted in these books).

In the study of lightlike hypersurfaces, the problem of 
construction of their normalizations and finding affine 
connections on such hypersurfaces arises naturally. 
This problem does not arise for the spacelike 
and timelike hypersurfaces s ince on them a family of normals is defined 
intrinsically in a first-order neighborhood: their 
normals are polar-conjugate of tangent hyperplanes $T_x (V), \;
x \in V$, with respect to the isotropic cones $C_x$ of 
the manifold $(M, g)$. For a lightlike hypersurface, 
a hyperplane $T_x (V)$ is tangent to the cone $C_x$. Hence 
a straight line orthogonal to $T_x (V)$ belongs to $T_x (V)$, 
and the family of these straight lines does not determine 
a normalization of a lightlike hypersurface $V$ and 
consequently an affine connection on $V$.

For a normalization of a lightlike hypersurface $V \subset 
(M, g)$ some authors (see  
   \cite{B72}, \cite{C65}, \cite{Ga67}, 
 \cite{K80}, \cite{SS38}) assign a field $N$ of isotropic 
directions not belonging to the tangent hyperplanes $T_x (V)$. 
Other authors (see, for example, the papers \cite{Be96a} and  
\cite{Be96b} and the book \cite{DB96}) assign a 
screen distribution $S$ on $V$ which belongs 
to the tangent bundle $T (V)$. Since an isotropic direction 
$N_x$ at a point $x \in V$ can be chosen being conjugate 
to a screen subspace $S_x$ with respect to 
the isotropic cone $C_x$, these two methods of 
normalization of a lightlike hypersurface $V \subset 
(M, g)$ are equivalent.

The important problem is to construct on a lightlike hypersurface 
$V \subset (M, g)$ a field of $N$ of isotropic directions 
or a screen distribution $S$ intrinsically connected with the 
geometry of $V$. Such a problem was open until now.

In this paper we present a few methods 
of construction of an invariant normalization 
on  a lightlike hypersurface $V$ of a  pseudo-Riemannian manifold 
$(M, g)$ of Lorentzian signature which is intrinsically connected 
with the geometry of $V$. In these constructions we 
use relative and absolute invariants defined by the first 
and second fundamental forms of $V$. The normalizations 
we have constructed are defined in a third-order neighborhood 
of a point $x$ of a lightlike hypersurface $V$. Each 
of the constructed normalizations 
induces an  affine connection  whose curvature tensor is 
expressed in terms of quantities connected with 
a fourth-order neighborhood of a point $x \in V$.

We describe briefly the contents of the paper. 
In Sections {\bf 1--3} we give the basic 
equations of the manifold 
$(M, g)$ of Lorentzian signature and construct on $(M, g)$ 
an isotropic frame bundle.  In Sections {\bf 4--5} 
we consider lightlike hypersurfaces $V$ on a manifold 
$(M, g)$, construct an isotropic frame bundle 
on them, and present the existence theorem 
for lightlike hypersurfaces.
 In Section {\bf 6} we study the fibration of 
isotropic geodesics on a lightlike hypersurface $V$, 
singular points and singular submanifolds of $V$.
 In Section {\bf 7} we find conditions defining 
invariant normalizations and affine connections on $V$.

Using the first 
and second fundamental forms of $V$, in Section {\bf 8} 
we construct on $V$ a series of relative and absolute 
invariants connected with 
a second-order neighborhood of a point $x \in V$. 
In Section {\bf 9} we consider the isotropic sectional 
curvature defined by Harris in \cite{H82}; see also  
\cite{BEE96}). 

Sections {\bf 10--11} are devoted to the construction 
of  invariant normalizations intrinsically connected with the 
geometry of a lightlike hypersurface $V$. As we have indicated 
earlier, these normalizations are constructed by means of 
the invariants that were found in  Section {\bf 8}, and they 
are defined in a third-order neighborhood of a point $x \in V$.

In the following two sections we address the problem of 
construction of an invariant normalization and 
an affine connection on lightlike hypersurfaces 
of some special classes: totally geodesic, totally umbilical, 
and belonging to a pseudo-Riemannian manifold of 
constant curvature. In these sections we clarify the 
role of the isotropic sectional curvature in 
the geometry of such hypersurfaces.

Note that in the paper [Be 96] and in Chapter 4 of the book [DB 96] 
for a lightlike hypersurface of  a pseudo-Riemanninan 
manifold $(M, g)$ (in particular, in a semi-Euclidean space 
 ${\bf R}^n_q$), a rigging (it is called a canonical screen distribution 
for ${\bf R}^n_q)$ and an induced affine connection 
 have been constructed. However, the authors did not give the
proof of independence of the constructed distribution and connection 
relative to a 
choice of a coordinate system in $(M, g)$ (in ${\bf R}^n_1$),
that is, they did not prove that 
these distribution and connection are intrinsically connected with the 
geometry of $V$.

Finally, in Section {\bf 14}, we consider a construction of an 
intrinsic normalization and an intrinsic affine connection 
on  lightlike hypersurfaces $V$ of a four-dimensional 
manifold $(M, g)$ of Lorentzian signature. We prove 
that in general, one can construct three normalizations and 
affine connections intrinsically connected with the 
geometry of $V$. Since a four-dimensional 
manifold $(M, g)$ of Lorentzian signature is directly 
connected with general relativity, the invariant 
normalizations we have constructed can have a physical meaning. 
In order to clarify the physical meaning, an assistance 
from physicists is needed.

In our study of lightlike hypersurfaces $V \subset (M, g)$ 
we use the method of moving frames and exterior 
differential forms of \'{E}. Cartan (see, for example, 
\cite{BCGGG91}, \cite{Ca67}, and \cite{AG93}). This allows us 
 to shorten computations and clarify a geometric 
meaning of constructed objects which is much more difficult 
in other methods.

The contents of this paper is directly connected 
with  our papers \cite{AG97}, \cite{AG98a}, \cite{AG98b}, 
\cite{AG98c},   and \cite{AGs} where we studied lightlike 
hypersurfaces in a pseudoconformal space, 
 the de Sitter space and on a manifold endowed with a conformal 
structure.

\setcounter{equation}{0}

{\bf 1. Pseudo-Riemannian manifolds of Lorentzian 
signature.}
Consider an $n$-dimensional 
 pseudo-Riemannian manifold $(M, g)$ of Lorentzian signature, 
where $M$ is a differentiable manifold of dimension 
$n,\; \dim M = n,$ and $g$ is a metric differential 
quadratic form of signature 
$(n-1, 1), \;\mbox{{\rm sign}} \; g = (n - 1, 1)$ 
(for definition see \cite{ON83}). 

A local frame associated with $(M, g)$ consists 
of a point $x \in M$ and $n$ vectors $e_i \in T_x (M), \;
i = 1, \ldots, n$, 
where $T_x (M)$ is a pseudo-Euclidean space tangent to 
the manifold $M$ at a point $x$. 

 For any two vectors $\xi, \;\eta \subset T_x (M), \;
\xi = \xi^i e_i, \; \eta = \eta^i e_i$, the quadratic form $g$ 
defines the scalar product  
 \begin{equation}\label{eq:1}
(\xi, \eta) = g (\xi, \eta) = g_{ij} \xi^i \eta^j,
\end{equation}
where  $g_{ij} = (e_i, e_j)$.

The equation 
 \begin{equation}\label{eq:2}
g (\xi, \xi) = 0 
\end{equation}
determines 
an  isotropic cone $C_x \subset T_x (M)$ at $x \in M$. 
The cone $C_x$ is real, and it carries 
rectilinear generators.

The equations of infinitesimal displacement of this 
frame have the form 
 \begin{equation}\label{eq:3}
dx = \omega^i e_i, \;\; de_i = \omega_i^j e_j, 
\end{equation}
where $\omega^i$ are basis forms of this manifold 
and $\omega^i_j$ are the  forms of the Levi-Civita 
connection. 

From (3) it follows that for a vector $\xi = \xi^i e_i$ we have
$$
d\xi = (d\xi^i + \xi^j \omega_j^i) e_i.
$$
The quantities 
$$
\nabla \xi^i = d\xi^i + \xi^j \omega_j^i
$$
are covariant derivatives of the coordinates of 
the vector $\xi$ in the Levi-Civita connection. The conditions of 
parallel displacement of the vector 
$\xi$ have the form $\nabla \xi^i = 0$. Since the scalar product 
remains unchanged under parallel displacement, we have 
$d (\xi, \eta) = 0$. It follows that in the Levi-Civita connection, the  metric tensor $g_{ij}$ satisfy 
 the following differential equations:
 \begin{equation}\label{eq:4}
\nabla g_{ij} = 
dg_{ij} - g_{ik} \omega_j^k - g_{kj} \omega_i^k = 0. 
\end{equation}
Equations (4) mean that the metric tensor is covariantly constant 
with respect to the Levi-Civita connection:

Note that the components $g_{ij}$ and the 1-forms $\omega^i$ 
are defined in a first-order differential neighborhood 
of a point $x \in (M, g)$, and the 1-forms $\omega_j^i$ 
are defined in its second-order neighborhood.


{\bf 2. The structure equations}.
 The forms  $\omega^i$ and $\omega^i_j$ are the  forms of the 
Levi-Civita connection. They 
 satisfy the following structure equations:
\begin{equation}\label{eq:5}
d\omega^i = \omega^j \wedge \omega^i_j, 
\;\; d\omega^i_j = \omega^k_j \wedge \omega^i_k +R^i_{jkl} 
\omega^k \wedge \omega^l,
\end{equation}
where   $i, j, k, l = 1,  \ldots, n$, and $R^i_{jkl}$ 
is the curvature tensor of the manifold $(M, g)$.  
The curvature tensor is defined in a third-order differential 
neighborhood of a point $x \in (M, g)$.

Consider the tensor 
 \begin{equation}\label{eq:6}
R_{ijkl} = g_{im} R^m_{jkl}.
\end{equation}
This tensor satisfies the following equations:
\begin{equation}\label{eq:7}
\renewcommand{\arraystretch}{1.3}
\left\{
\begin{array}{ll}
R_{ijkl} = -R_{jikl} = -R_{ijlk}, \\
R_{ijkl} = R_{klij}, \\
R_{ijkl} + R_{iklj} + R_{iljk} = 0.
\end{array}
\right.
\renewcommand{\arraystretch}{1}
\end{equation}

If the curvature tensor vanishes, $R^i_{jkl} = 0$, then 
$(M, g)$ is a pseudo-Euclidean space $R^n_1$ of signature 
$(n - 1, 1)$ (for $n = 4$, it is a Minkowski space), and equations (3) are completely integrable for such a space. 
 
If the curvature tensor does not vanish, $R^i_{jkl} \neq 0$, then 
equations (3) are integrable along a curve $x = x (t) \subset M$. 
A solution of these equations defines a {\em development} of this 
line and the frame bundle along the curve onto the tangent 
pseudo-Euclidean space $(R^n_1)_x$ at the point $x \in M$.

{\bf 3. An isotropic frame bundle on \protect\boldmath 
\((M, g)\). \protect\unboldmath} 
Let $C_x$ be an isotropic cone, let $\eta$ be an isotropic 
hyperplane, and let $e_1$ be an isotropic vector along which 
the hyperplane $\eta$ is tangent to the cone $C_x$. 
Let further the vectors $e_a \in \eta, \; a = 2, \ldots , 
n - 1$, be spacelike vectors, and let $e_n$ be an isotropic 
(normalizing) vector not belonging to $\eta$ 
and conjugate to the vector $e_a$. Suppose that 
$\zeta$ is a hyperplane tangent to $C_x$ along $e_n$. 
Then the $(n-2)$-dimensional subspace 
$S_x = \eta \cap \zeta = e_2 \wedge \ldots \wedge 
e_{n-1}$ is called a {\em screen subspace}. 


In the isotropic frame described above the matrix 
of the metric tensor $g$ has the form 
 \begin{equation}\label{eq:8}
 (g_{ij}) = \pmatrix{0 & 0& -1 \cr
           0 & g_{ab}& 0 \cr
           -1 & 0& 0 \cr}, \; a, b = 2, \ldots , n-1.
\end{equation}
Here $a, b = 2, \ldots, n - 1, \; g_{1n} = (e_1, e_n) = -1$ 
is a normalizing 
condition,  $\det (g_{ab}) \neq 0, \; \mbox{{\rm rank}} \; (g_{ab}) = n - 2$, and $g_{ab} \xi^a \xi^b > 0$.


It follows from equations (4) and (8) that 
\begin{equation}\label{eq:9}
g =  g_{ab} \xi^a \xi^b - 2 \xi^1 \xi^n,
\end{equation}
 \begin{equation}\label{eq:10}
\renewcommand{\arraystretch}{1.3}
\left\{
\begin{array}{ll}
\omega_1^n = \omega^1_n = 0, & \omega_1^1 + \omega_n^n = 0, \\ 
\omega_a^n = g_{ab} \omega_1^b, &\omega_a^1  = g_{ab} \omega_n^b, 
  \\
dg_{ab} - g_{ac} \omega_b^c - g_{cb} \omega_a^c = 0. &
\end{array}
\right.
\renewcommand{\arraystretch}{1}
\end{equation}

{\bf 4. Lightlike hypersurfaces.}
Suppose that $V \subset (M, g), \; \dim \; V = n - 1$,
 is a lightlike hypersurface 
on the manifold $(M, g)$, and $x \in V$ is a point of 
$V$. Then  the tangent hyperplane $\eta = T_x (V)$ is 
isotropic, i.e., it is tangent to the cone $C_x$. 
Let $e_1$ be an isotropic vector in $\eta$ which together with 
vectors $e_a, \;a = 2, \ldots , n - 1$, 
form a basis of the subspace $\eta$. Finally suppose 
that $e_n \notin \eta$ is also an isotropic vector 
(see Section {\bf 3}). 

Then the equation of $V$ is
\begin{equation}\label{eq:11}
\omega^n = 0.
\end{equation}
On the hypersurface $V$ we have 
\begin{equation}\label{eq:12}
g = g_{ab} \xi^a \xi^b, \;\; 
\mbox{{\rm rank}}  \;\; g = n - 2.
\end{equation}
This form is called the {\em first fundamental form} 
of $V$, and the equations $\omega^a = 0$ define 
{\em isotropic lines} on $V$. 

Consider a first-order frame bundle associated with a 
lightlike hypersurface $V \subset (M, g)$. Since by (3) 
and (11) we have 
\begin{equation}\label{eq:13}
dx = \omega^1 e_1 + \omega^a e_a,
\end{equation}
the forms $\omega^1$ and $\omega^a$ are basis forms on the 
hypersurface $V$. If we fix a point $x \in V$, we obtain that 
$\omega^1 =\omega^a = 0$. As a result, equations (3) 
take the form
\begin{equation}\label{eq:14}
\renewcommand{\arraystretch}{1.3}
\left\{
\begin{array}{lll}
\delta e_1 = &\!\!\!\!\pi_1^1 e_1, &\!\!\!\!\!\\
\delta e_a = &\!\!\!\! \pi_a^1 e_1 & \!\!\!\!\!
+ \pi_a^b e_b,\\
\delta e_n = &  
&\!\!\!\!\!  + \pi_n^a e_a - \pi_1^1 e_n,
\end{array}
\right.
\renewcommand{\arraystretch}{1}
\end{equation}
where  $\delta = d|_{\omega^1 = \omega^a = 0}$ 
is the symbol of differentiation with respect to fiber 
parameters and $\pi_\eta^\xi = \omega_\eta^\xi (\delta) 
= \omega_\eta^\xi|_{\omega^1 = \omega^a = 0}$. 

By (10), we find that 
\begin{equation}\label{eq:15}
\pi_n^a = g^{ab} \pi_b^1.
\end{equation}
Thus the forms 
$\pi_1^1, \pi_b^a$, and $\pi_a^1$  are independent fiber 
forms.  These forms are invariant 
forms of the group of admissible transformations 
of first-order frames whose dimension is $1 + (n - 2) 
+ (n - 2)^2 = n -1 + (n-2)^2$.

Among the fiber forms the forms $\pi_1^a$ play a special 
role. They define a displacement of a screen distribution $S_x$ in the tangent hyperplane $T_x (V)$ of a lightlike 
hypersurface $V$. By (15) 
there is a bijective correspondence between 
the screen subspaces $S_x$ and the normalizing isotropic 
straight lines $x e_n = N_x$.

Taking exterior derivatives of equation (10), we arrive at the 
exterior quadratic equation 
\begin{equation}\label{eq:16}
\omega^a \wedge \omega_a^n = 0. 
\end{equation}
Applying Cartan's lemma to this equation, we find that 
\begin{equation}\label{eq:17}
\omega_a^n = \lambda_{ab} \omega^a, \;\;\lambda_{ab} = \lambda_{ba}.
\end{equation}
The tensor $\lambda_{ab}$ forms the 
{\em second fundamental tensor}  
of the hypersurface  $V$, and its {\em second fundamental form}
 is
\begin{equation}\label{eq:18}
\varphi = \lambda_{ab} \omega^a \omega^b.
\end{equation}

Equations (10) imply that 
\begin{equation}\label{eq:19}
 \omega_1^a = \lambda_b^a \omega^b,
\end{equation}
where $\lambda_b^a = g^{ac} \lambda_{cb}$ is the 
{\em Burali-Forti affinor} of $V$ (see \cite{Bu12}). 
Note that the authors 
of \cite{DB96} called $\lambda_b^a$ the shape operator 
(see \cite{DB96}, pp. 85, 154, and 160).  

Equations (3), (10), and (11) imply that 
\begin{equation}\label{eq:20}
de_1 =  \omega^1_1 e_1 + \omega_1^a e_a.
\end{equation}
The point $x$ and the vector $e_1$ define an isotropic 
direction $x e_1$ on the hypersurface $V$. 
 By (19) the system of equations $\omega^a = 0$ defines an isotropic 
fibration ${\cal F}$ on $V$ and $V = M^{n-2} \times l$, where 
$l$ is a straight line whose image is 
an isotropic geodesic $x e_1$ on the 
manifold $(M, g)$,  $f (l) = x e_1$  (see \cite{AGs}). 

{\bf 5. The existence theorem.} 
Applying the Cartan test (see \cite{BCGGG91}) to the system 
of equations (11), (16), and (17) in the same way as in 
\cite{AGs}, we arrive at the following theorem.

\begin{theorem} Lightlike hypersurfaces on a 
manifold $(M, g)$ exist, and the solution of a system 
defining such hypersurfaces depends on one function 
of $n - 2$ variables.
\end{theorem}

{\sf Proof}. The proof of Theorem 1 coincides with 
the proof of the existence theorem for lightlike hypersurfaces 
$V$ on a manifold $(M, c)$ endowed with a conformal structure 
of Lorentzian signature given in \cite{AGs}. \rule{3mm}{3mm}

{\bf 6. Isotropic geodesics on \protect\boldmath \( V \subset (M, g). \) \protect\unboldmath} 
It follows from (12) and (18) that integral curves $\gamma$ 
of the vector field $e_1$ defined by the equations 
$\omega^a = 0$ are isotropic and asymptotic on $V$. 
These curves form a foliation ${\cal F}$ on $V$. 

\begin{theorem}
Isotropic lines $\gamma$ of a lightlike hypersurface 
$V$ are geodesic lines of the manifold $(M, g)$. 
\end{theorem}

{\sf Proof}. In fact, the equations of geodesic lines on 
a Riemannian manifold have the form
\begin{equation}\label{eq:21}
d \omega^i + \omega^j \omega_j^i = \alpha \omega^i,
\end{equation}
where $\alpha$ is an 1-form. 
For $i = a$, these equations becomes 
$$
d \omega^a + \omega^1 \omega_1^a + \omega^b \omega_b^a 
= \alpha \omega^a.
$$
It follows from (19) that 
 for $\omega^a = 0$, equations (21) are satisfied identically. 
\rule{3mm}{3mm}

Note that the isotropic geodesics on pseudo-Riemannian manifolds were
considered in \cite{AG97} (see also \cite{AG96}), where, in particular,
their invariance under conformal transformations of a pseudo-Riemannian metric has been proved. 

Theorem 2 implies that {\em the foliation ${\cal F}$ 
is also a geodesic foliation on $V$.} 

Under the development of the manifold $(M, g)$ onto 
the tangent pseudo-Euclidean space $(R^n_1)_x = T_x(M)$, 
to the isotropic geodesic $x e_1$ there corresponds 
the straight line $l$. 
Consider a point $y = x + s e_1$ on the straight line 
$l$. From equations (20) it follows that 
$$
d y = (ds + s \omega^1_1 + \omega^1) e_1 + (\omega^a 
+ s\omega^1_a) e_a.
$$
But by (19), we have 
$$
\omega^a + s\omega^1_a = (\delta_b^a + s \lambda_b^a) \omega^b.
$$
This allows us to write the equation for $dy$ in the form
\begin{equation}\label{eq:22}
d y = (ds + s \omega^1_1 + \omega^1) e_1 + (\delta_b^a 
+ s \lambda_b^a) \omega^b e_a.
\end{equation}

The matrix $(J_b^a) = (\lambda_b^a + s \delta_b^a)$ 
is the Jacobi matrix of the mapping \newline 
$f: M^{n-2} \times l \rightarrow V \subset (M, g)$, and 
its determinant,
$$
J = \det (\lambda_b^a + s \delta_b^a)
$$
is the Jacobian of this mapping. 

Since the affinor $\lambda_b^a = g^{ac} \lambda_{cb}$ 
is symmetric, its characteristic equation 
\begin{equation}\label{eq:23}
\det (\lambda_b^a - \lambda \delta_b^a) = 0
\end{equation}
has $n - 2$ real roots $\lambda_a$ if each of them is counted as 
many times as its multiplicity. This implies the following 
theorem.

\begin{theorem} Any isotropic geodesic $l$ of 
a lightlike hypersurface  $V$ of a manifold $(M, g)$ 
carries $n - 2$ real singular points 
if each of them is counted as many times 
as its multiplicity.
\end{theorem}

{\sf Proof.} Consider the development $\widetilde{V}$  of the 
hypersurface $V$ 
onto the tangent  space $(R^n_1)_x = T_x(M)$. 
The tangent subspace $T_y (\widetilde{V})$ to 
  $\widetilde{V}$ at a point $y$ 
is a subspace of the space $T_x (M)$. By (22), 
this subspace is determined by the point $y$ and 
the vectors $e_1$ and $f_b = (\lambda_b^a + s \delta_b^a) e_a$. If the Jacobian $J$ is different from 0, then these vectors 
are linearly independent and 
determine the $(n-1)$-dimensional tangent subspace 
$T_y (V)$. In this case the point $y$ is a regular 
point of the hypersurface  $\widetilde{V}$, 
and to such a point on $\widetilde{V}$  there corresponds 
 a regular point of $V \subset (M, g)$. 
If at a point $y \in x e_1$ the Jacobian $J$ is equal to 0, 
then at this point $\dim T_y (\widetilde{V}) < n - 1$, 
and this point is a singular point of $\widetilde{V}$. 
To such a point on $\widetilde{V}$ there corresponds 
a singular point of the hypersurface $V \subset (M, g)$.

 Singular points are 
defined by the equation 
\begin{equation}\label{eq:24}
\det (\lambda_b^a + s \delta_b^a) = 0
\end{equation}
Comparing equations (23) and (24), 
we find the coordinates $s_a$  of these 
singular points: $s_a = - \frac{1}{\lambda_a}$. 
Thus the singular points of the straight line $l$ are 
\begin{equation}\label{eq:25}
F_a = x - \frac{1}{\lambda_a} e_1. \rule{3mm}{3mm}
\end{equation}

Note that if $\lambda_a = 0$, then $F_a$ is the point 
at infinity.
It is obvious that the point $x$ is a regular point 
of the straight line $l$.

To an eigenvalue $\lambda_a$ of the affinor 
$(\lambda_b^a)$ there corresponds an invariant two-dimensional 
eigenplane passing through the vector $e_1$. The eigenplanes 
corresponding to distinct eigenvalues $\lambda_a$ and $\lambda_b 
\neq \lambda_a$ are orthogonal with respect to the scalar 
product $(\xi, \eta) = g_{ab} \xi^a \eta^b$. 

If $\lambda_a$ is a simple root of equation (23), then 
the focus $F_a$ describes a {\em lightlike focal 
submanifold} $(F_a), 
\dim \; (F_a) = n - 2$, carrying an $(n-3)$-parameter 
family of isotropic lines. The eigenplane corresponding to 
such a root $\lambda_a$ is osculating for these lines. 

In the paper \cite{AG98a}, for a lightlike hypersurface 
of a pseudo-Riemannian de Sitter space 
we investigated the structure of such singular points 
 and the structure of $V$ itself taking into account 
multiplicities of singular points.  Many of the results of 
\cite{AG98a} are still valid for a lightlike hypersurface  $V 
\subset (M, g)$.

{\bf 7. An affine connection on $V \subset (M, g)$.} 
From equations (5) it follows that the basis forms 
$\omega^1$ and $\omega^a$ of the hypersurface $V$ satisfy
the following structure equations:
\begin{equation}\label{eq:26}
\renewcommand{\arraystretch}{1.3}
\left\{
\begin{array}{ll}
d \omega^1 = \omega^1 \wedge \omega^1_1 + 
\omega^a \wedge \omega^1_a, \\
d \omega^a = \omega^1 \wedge \omega^a_1 + 
\omega^b \wedge \omega^a_b.
\end{array}
\right.
\renewcommand{\arraystretch}{1}
\end{equation}
Thus the 1-form
$$
\omega = \pmatrix{\omega^1_1 &  \omega_a^1 \cr 
\omega^a_1 & \omega_b^a}
$$
defines an affine structure on $V$.  
 To define an 
 affine connection, the form $\omega$ must satisfy the 
structure equation
\begin{equation}\label{eq:27}
d\omega +  \omega \wedge \omega = \Omega,
\end{equation}
where $\Omega$ is the curvature 2-form of this connection 
which is a linear combination of exterior products 
of the basis forms $\omega^1$ and $\omega^a$ 
(see, for example, \cite{KN63}, Ch. III).

Taking the exterior derivative of the form $\omega$ 
componentwise and applying equations (5), (10), and (11), 
we find that 
\begin{equation}\label{eq:28}
\renewcommand{\arraystretch}{1.3}
\left\{
\begin{array}{ll}
d \omega^1_1 + \omega_a^1 \wedge \omega^a_1 = R^1_{1kl} \omega^k \wedge \omega^l, \\ 
d \omega_a^1 + \omega^1_1 \wedge \omega^1_a + \omega_b^1 \wedge 
\omega_a^b = R^1_{akl} \omega^k \wedge \omega^l, \\ 
d \omega^a_1 + \omega^a_1 \wedge \omega^1_1 + \omega_b^a \wedge 
\omega_1^b =  R^a_{1kl} \omega^k \wedge \omega^l, \\ 
d \omega^a_b + \omega^a_1 \wedge \omega^1_b 
+ \omega_c^a \wedge \omega_b^c  
= \omega_b^n \wedge \omega_n^a  
+ R^a_{bkl} \omega^k \wedge \omega^l. 
\end{array}
\right.
\renewcommand{\arraystretch}{1}
\end{equation}

Equations (28) and (17) show that  conditions (27) 
are satisfied if and only if the 1-form  
 $\omega_a^1$, and by (10) the form 
$\omega_n^a$ as well, are expressed in terms of the basis forms 
of the hypersurface $V$: 
\begin{equation}\label{eq:29}
\omega^1_a = \nu_a  \omega^1 + \nu_{ab}  \omega^b, \;\; 
 \omega_n^a = g^{ab}  \omega^1_b. 
\end{equation}

It follows from (3) that the vectors $e_a$ and $e_{n}$ 
satisfy the differential equations 
 \begin{equation}\label{eq:30}
\renewcommand{\arraystretch}{1.3}
\left\{
\begin{array}{ll}
de_a = \omega_a^1 e_1 + & \!\!\!\! 
\omega^b_a e_b + \omega_a^n e_n,  \\ 
de_n =   &\!\!\!\!  \omega^a_n e_a - \omega_1^1 e_n.
\end{array}
\right.
\renewcommand{\arraystretch}{1}
\end{equation}
For $\omega^1 = \omega^a = 0$, equations (30) 
take the form 
 \begin{equation}\label{eq:31}
de_a = \omega^b_a e_b,  \;\; de_n = - \omega_1^1 e_n.
\end{equation}
This means that conditions (29) are satisfied if and 
only if the screen distribution $S = \cup_{x \in V} S_x$, 
or equivalently the field of normalizing isotropic straight 
lines $N = \cup_{x \in V} x e_1$, 
are defined invariantly. Note the summation in these two 
expressions are carried over the regular points $x \in V$. 

Hence  {\em an affine connection on $V$ is defined 
if and only if on $V$ there is given an invariant screen 
distribution $S$ (or a field of normalizing isotropic straight 
lines $N$).} This result is well-known and was discussed 
in many papers. Note that Bonnor \cite{B72}, Cagnac \cite{C65}, 
Galstyan \cite{Ga67}, Katsuno \cite{K80},   Lemmer \cite{L65}  
(see also \cite{SS38}) 
constructed a field of isotropic normalizing vectors while 
Duggal and Bejancu in their book \cite{DB96} considered a screen 
distribution.

However, in all papers on this subject known to the authors, 
the problem of construction of a screen 
distribution $S$ or a field of normalizing isotropic straight 
lines $N$ that are {\em intrinsically connected} 
with the geometry 
of a lightlike hypersurface $V \in (M, g)$ was not considered. 
In what follows we present a few solutions of this problem.

{\bf 8. Invariants of a lightlike hypersurface.} 
A lightlike hypersurface $V \subset (M, g)$ in an isotropic 
first order frame is determined by equation (11) whose 
prolongation gives equation (17). 

Exterior differentiation of equations (17) by means of structure 
equations (5) and equations (10) leads to 
the following exterior quadratic equations:
$$
[\nabla \lambda_{ab} - \lambda_{ab} \omega_1^1 
+ (\lambda_{ac} g^{ce} \lambda_{eb} 
+ 2 R^n_{ab1}) \omega^1 + R^n_{abc} \omega^c] \wedge 
 \omega^b = 0, 
$$
where $\nabla \lambda_{ab}= d \lambda_{ab} - 
\lambda_{ac}  \omega_b^c - \lambda_{cb} \omega_a^c$. 
Applying Cartan's lemma to the last equation, we find that  
\begin{equation}\label{eq:32}
\nabla \lambda_{ab} - \lambda_{ab} \omega_1^1 + (\lambda _{ac} 
g^{ce} \lambda_{eb} + 2 R^n_{ab1}) \omega^1 + R^n_{abc} \omega^c 
 = \mu_{abc} \omega^c. 
\end{equation}
Here the quantities $\mu_{abc}$ are symmetric with respect 
to all indices.

The quantities $R^n_{ab1}$ are symmetric with 
respect to the indices $a$ and $b$  since by (6) and (7) we have
$$
R^n_{ab1} = - R_{1ab1} = -  R_{b11a} = - R_{1ba1} = R^n_{ba1}.
$$
Now if we alternate equations (32) with respect 
to the indices $a$ and $b$, then we find that 
$$
R^n_{[ab]c} = 0.
$$
This implies
$$
R^n_{abc} = R^n_{bac}.
$$
But since by (7) we have
$$
R^n_{abc} = - R^n_{acb},
$$
we find that 
$$
R^n_{abc}  = - R^n_{acb} = - R^n_{cab} = R^n_{cba}  
= R^n_{bca} = - R^n_{bac} = - R^n_{abc}.
$$
It follows that 
\begin{equation}\label{eq:33}
R^n_{abc} = 0.
\end{equation}

Hence on a lightlike hypersurface $V \subset (M, g)$ 
conditions (33) are satisfied. As a result, equations (32) 
take the form
\begin{equation}\label{eq:34}
\nabla \lambda_{ab} - \lambda_{ab} \omega_1^1 + (\lambda _{ac} 
g^{ce} \lambda_{eb} + 2 R^n_{ab1}) \omega^1 
 = \mu_{abc} \omega^c. 
\end{equation}

For a fixed point $x \in V$ (i.e., for 
$\omega^1 = \omega^a = 0$), we find from (34) that 
\begin{equation}\label{eq:35}
\nabla_\delta \lambda_{ab} = \lambda_{ab} \pi_1^1,
\end{equation}
where  
$$
\nabla_\delta \lambda_{ab} = \delta \lambda_{ab} 
- \lambda_{ac} \pi_b^c - \lambda_{cb} \pi_a^c.
$$

Equations (35) prove that the quantities $\lambda_{ab}$ 
 form a relative (0, 2)-tensor of weight 1. This tensor 
is the {\em second fundamental tensor} of 
the hypersurface $V$. It is defined in a second-order 
neighborhood of a point $x \in V$.

It follows from (10) and (35) that for a fixed point $x \in V$ 
the affinor $\lambda_b^a$ satisfies the equations
\begin{equation}\label{eq:36}
\nabla_\delta \lambda_b^a = \lambda_b^a \pi_1^1.
\end{equation}
Hence  it is also of weight 1.

Consider characteristic equation (23) of the affinor 
$\lambda_b^a$. We write it in the expanded form
\begin{equation}\label{eq:37}
\lambda^{n-2} - I_1 \lambda^{n-3} + \ldots 
+ (-1)^{n-2} I_{n-2} = 0.
\end{equation}
The coefficients of this equation 
are relative invariants of weights equal 
to their labels. These invariants are the sums of 
the diagonal minors of corresponding orders 
of the matrix $(\lambda_b^a)$:
\begin{equation}\label{eq:38}
I_1 = \lambda_a^a, \;\; I_2 = \lambda_{[a}^b \lambda_{b]}^a, 
\ldots , \; I_{n-2} = \det (\lambda_b^a).
\end{equation}
These coefficients form a complete system of relative 
invariants of the affinor $\lambda_b^a$.  
We can get  another complete system of  relative 
invariants of the affinor $\lambda_b^a$ if we consider 
the following contractions:
\begin{equation}\label{eq:39}
\widetilde{I}_1 = I_1 = \lambda_a^a, \;\; 
\widetilde{I}_2 = \lambda_a^b \lambda_b^a, \ldots , 
\widetilde{I}_{n-2} =  \lambda_{a_1}^{a_{n-2}}   \lambda_{a_2}^{a_1} \ldots   \lambda_{a_{n-2}}^{a_{n-3}}.
\end{equation}

Moreover, the roots $\lambda_a, \; a = 2, \ldots , n - 1$, 
 of characteristic equation (37) also form 
a complete system of  
invariants of weights 1 of the affinor $\lambda_b^a$.  

We can find invariants of weights 1 from 
nonvanishing invariants (38) and (39) 
if we take from them the root of degree equal to their labels: 
the quantities $|I_p|^\frac{1}{p}$ and  
$|\widetilde{I}_p|^\frac{1}{p}$ are invariants of weight 1.

Equations (36) imply that for a fixed point $x \in V$, 
each relative invariant $I$ of weight 1 satisfies 
the differential equation
\begin{equation}\label{eq:40}
\delta I = I \pi_1^1.
\end{equation}
Any nonvanishing relative invariant $I$ of weight 1 
allows us to normalize the isotropic vector $e_1$ by setting 
$\widetilde{e}_1 = \displaystyle 
\frac{1}{I} e_1$, and the new vector $\widetilde{e}_1$ 
is invariant. In fact, it follows from 
(3) and (10) that for a fixed point $x \in V$ we have
$$
\delta e_1 = \pi_1^1 e_1.
$$
This and equation (40) imply that $\delta \widetilde{e}_1 = 0$, 
and thus the vector $\widetilde{e}_1$ does not depend on 
a choice of normalizing parameter on an isotropic 
geodesic $x e_1$.

Absolute invariants of a hypersurface $V$ can be constructed by 
taking ratios of two nonvanishing relative invariants of 
the same weight. For a fixed point $x \in V$, 
an absolute invariant $J$ satisfies the equation
\begin{equation}\label{eq:41}
\delta J = 0.
\end{equation}

Since  the affinor $\lambda_b^a$ 
is defined in a second-order neighborhood of a point $x \in V$, 
it follows that all absolute 
and relative  invariants of a hypersurface $V$ constructed by 
means of  $\lambda_b^a$  are defined also 
in a second-order neighborhood of  $x \in V$,

{\bf 9. Isotropic sectional curvature of a lightlike 
hypersurface.} 
Harris    introduced 
the notion of isotropic sectional curvature 
of an isotropic 2-plane $\sigma$ of a pseudo-Riemannian 
manifold $(M, g)$ (see  \cite{H82}; see 
also  the book \cite{BEE96}, Appendix A, p. 571). 
 If $N$ is a isotropic 
nonzero element of a one-dimensional space 
of isotropic vectors belonging to $\sigma$, and 
$P$ is an arbitrary (nonzero) nonisotropic vector from $\sigma$, 
then the isotropic sectional curvature $K_N (\sigma)$ is defined 
as 
\begin{equation}\label{eq:42}
K_N (\sigma) = \frac{(R (P, N) N, P)}{(P, P)}.
\end{equation}
This expression does not depend on a vector $P \subset \sigma$ 
but depends quadratically on an isotropic  vector $N$.

Denote by $n^i$ coordinates of an isotropic  vector $N$ and 
by $p^i$ coordinates of a vector $P$. Then for the standard 
coordinate representation of the curvature tensor 
(see (5) and (6)) the nominator of (42) can be written as 
$$
(R (P, N) N, P) = R_{ijkl} n^i p^j p^k n^l,
$$
and its denominator is $(P, P) = g_{ij} p^i p^j$. 

Let $V$ be a lightlike hypersurface of a pseudo-Riemannian 
manifold $(M, g)$ of Lorentzian signature, and let $T_x (V)$ 
be its tangent hyperplane. In the isotropic frame 
considered in section {\bf 4}, the vector $e_1$ is 
isotropic, and this vector and a vector $P = p^1 e_1 + p^a e_a$ 
determine an isotropic 2-plane $\sigma = e_1 \wedge P$. 
For this 2-plane the isotropic sectional curvature has 
the following expression:
\begin{equation}\label{eq:43}
K_N (\sigma) = \frac{R_{1ab1} p^a p^b}{g_{ab} p^a p^b}.
\end{equation}

A  lightlike hypersurface $V \subset (M, g)$ is called a 
{\em hypersurface of null isotropic sectional curvature} if 
for all its tangent two-dimensional isotropic planes $\sigma$, 
their isotropic sectional curvatures vanish. 

Consider equation (32) for the second fundamental tensor $\lambda_{ab}$ of a lightlike hypersurface $V \subset 
(M, g)$. This equation contains the components $R^n_{ab1}$ of 
the curvature tensor of the manifold $(M, g)$. But by 
(7) we have 
\begin{equation}\label{eq:44}
R_{1ab1} = - R^n_{ab1}.
\end{equation}

Now we will prove the following theorem. 

\begin{theorem}
The isotropic sectional curvature of a lightlike hypersurface 
$V \subset (M, g)$ vanishes if and only if the derivative of the 
second fundamental tensor of $V$ along the field of isotropic 
directions on $V$ is expressed in terms of elements 
of a second-order neighborhood.
\end{theorem}

{\sf Proof.} The field of isotropic directions on $V$ is 
defined by the equations $\omega^a = 0$. It follows from 
equation (34) that the derivative of the tensor $\lambda_{ab}$ along an isotropic direction on $V$ is determined by 
the formula 
\begin{equation}\label{eq:45}
(\nabla \lambda_{ab} - \lambda_{ab} \omega_1^1)_{,1} 
= - \lambda _{ac} g^{ce} \lambda_{eb} - 2 R^n_{ab1}. 
\end{equation}
In the right-hand side of this equation the first 
term is defined in a second-order neighborhood 
of a point $x \in V$, and the second term in its 
third-order neighborhood. By (43) and (44), the second 
term vanishes if and only if a hypersurface $V$ 
has its isotropic sectional curvature equal to 0. 
\rule{3mm}{3mm}

It follows from Theorem 4 that the derivatives of 
all the invariants of a lightlike hypersurface 
with the  vanishing isotropic sectional curvature taking along a 
field of isotropic directions of $V$ are also defined in terms 
of second-order objects.

{\bf 10. Construction of a screen distribution by means of 
absolute invariants.} We will prove the following theorem.

\begin{theorem}
If $J = J (x)$ is an absolute invariant defined on 
a lightlike hypersurface $V \subset (M, g)$, and the level 
$(n-2)$-dimensional submanifold of $J (x)$ are transversal 
 to isotropic geodesics of $V$, then the distribution 
$S$ tangent to these level  submanifolds is an invariant screen 
distribution. If the invariant $J (x)$ is connected with 
the hypersurface $V$ intrinsically, then the same is true for 
a screen distribution $S$ generated by $J$. 
If the order of an invariant $J (x)$ is equal 
to $p$, then the normalization  is defined in a 
 neighborhood of a point $x \in V$ of order $p + 1$, 
and the curvature tensor in a 
 neighborhood of a point $x \in V$ of order $p + 2$. 
\end{theorem}

{\sf Proof.} By (41), the differential of the invariant $J$ 
has the form
\begin{equation}\label{eq:46}
d J = K \omega^1 + \widetilde{K}_a \omega^a,
\end{equation}
where $K \neq 0$. On a level  submanifold, 
$d J = 0$. It follows that 
\begin{equation}\label{eq:47}
 \omega^1 = K_a \omega^a,
\end{equation}
where $K_a = - \displaystyle\frac{\widetilde{K}_a}{K}$. Thus on a level 
surface we have 
$$
dx = \omega^a \widetilde{e}_a,
$$
where $\widetilde{e}_a = e_a + K_a e_1$. At a point $x \in V$, 
the vectors $\widetilde{e}_a$ determine an invariant screen 
subspace $S_x = \widetilde{e}_2 \wedge \widetilde{e}_3 \wedge 
\ldots \wedge \widetilde{e}_{n-1}$. The distribution 
$S = \cup_{x \in V} S_x$ is an invariant screen distribution 
generated by the invariant $J = J(x)$. If this invariant 
is intrinsically connected with the hypersurface $V$, then 
the same is true for the screen distribution $S$ generated by 
$J$. 

Let us make a reduction in the isotropic first-order 
frame bundle by superposing the vectors $e_a$ with 
the vectors $\widetilde{e}_a$. Then we have $K_a = 0$, and 
equation (47) takes the form
$$
 \omega^1 = 0.
$$
Since this equation determines a family of level  submanifolds of 
the invariant $J$, it must be completely integrable. Hence
$$
 d\omega^1 \wedge  \omega^1 = 0.
$$
By (5), the last equation can be written as
$$
 \omega^1 \wedge  \omega^a  \wedge  \omega^1_a = 0.
$$
This implies that 
\begin{equation}\label{eq:48}
 \omega_a^1 = \nu_a \omega^1 + \nu_{ab} \omega^a,
\end{equation}
where $\nu_{ab} = \nu_{ba}$. Equation (48) coincides with 
equation (29). However the condition $\nu_{ab} = \nu_{ba}$ 
shows that an affine connection generated by an 
absolute invariant $J$ is a connection of special type.
If an absolute invariant $J = J(x)$ is constructed by means 
of the affinor $\lambda_b^a$, then it is defined in a 
second-order neighborhood of a point $x \in V$, the 
quantities $K_a$ and $\widetilde{K}_a$ defining the screen 
distribution in a 
third-order neighborhood, and finally, the quantities 
$\nu$ and $\nu_a$ from equations (48) in a 
fourth-order neighborhood. 
Thus the curvature tensor of the 
affine connection generated by the absolute invariant $J$ is also 
defined in a fourth-order neighborhood of a point $x \in V$.
 \rule{3mm}{3mm}

{\bf 11. Construction of a screen distribution by means of 
relative invariants.} In a first-order frame bundle of a lightlike hypersurface $V$ constructed in Section {\bf 4}, we define a  screen subspace $S_x$ by vectors $c_a$:
$$
c_a = e_a + z_a e_1, \;\;\;\; a = 2, \ldots , n-1.
$$
This subspace is invariant if and only if 
\begin{equation}\label{eq:49}
\delta c_a = \sigma_a^b c_b,
\end{equation}
where as earlier, $\delta$ is the symbol of differentiation 
with respect to fiber parameters, and $\sigma_a^b$ are 
some 1-forms. 

Applying equations (3), (10), (11), and (17), we find that 
\begin{equation}\label{eq:50}
\delta c_a = (\nabla_\delta z_a + z_a \pi^1_1 + \pi_a^1) e_1 
+ \pi_a^b c_b.
\end{equation}
Comparing equations (50) and (49), we see that the screen 
subspace $S_x = [x, c_2, \ldots , c_{n-1}]$ is invariant if and 
only if the following conditions hold:
\begin{equation}\label{eq:51}
\nabla_\delta z_a + z_a \pi^1_1 + \pi_a^1 = 0.
\end{equation}
The coordinates of a {\em normalizing object} $z_a$ 
defining an invariant screen subspace $S_x$ 
must satisfy this equation.

Consider a nonvanishing relative invariant $I = I (x)$ of weight 
1 defined in a second-order neighborhood of a point $x \in V$. 
Equation (40) which this invariant satisfies can be written as 
$$
\delta \ln |I| = \pi_1^1.
$$
The last equation is equivalent to 
the equation 

\begin{equation}\label{eq:52}
d \ln |I| -  \omega_1^1 = - K \omega^1 - K_a \omega^a.
\end{equation}
The coefficients $K$ and $K_a$ in (52) are 
defined in a third-order neighborhood of a point $x \in V$. 

We will prove the following theorem.

\begin{theorem}
If the coefficient $K$ in 
equation $(52)$  is not  
 a root of characteristic equation $(37)$, then 
the  coefficients  $K_a$ in equation $(52)$ 
allows one to construct an object  defining 
an invariant normalization of a lightlike 
hypersurface $V \subset (M, g)$. This normalization is 
intrinsically connected with the geometry of $V$ and 
defined in a third-order neighborhood of 
a point $x \in V$.
\end{theorem}

{\sf Proof.} 
Taking exterior derivatives of equation (52), we find that 
 \begin{equation}\label{eq:53}
\renewcommand{\arraystretch}{1.3}
\begin{array}{ll}
(dK -  K \omega_1^1) \wedge  \omega^1 
+ (\nabla K_a + (\lambda_a^b - K \delta_a^b) 
\omega_b^1) \wedge \omega^a \\
+ K_b \lambda_a^b  \omega^1 \wedge \omega^a - 
R^1_{1kl}  \omega^k  \wedge \omega^l = 0, 
\end{array}
\renewcommand{\arraystretch}{1}
\end{equation}
where $\nabla K_a = dK_a - K_b \omega^b_a$. 
It follows from equation (53) that 
 \begin{equation}\label{eq:54}
\renewcommand{\arraystretch}{1.3}
\left\{
\begin{array}{ll}
dK -  K \omega_1^1 = M \omega^1 + M_a \omega^a,\\
\nabla K_a + (\lambda_a^b - K \delta_a^b) \omega_b^1 
= \widetilde{M}_a  \omega^1 + M_{ab}  \omega^b.
\end{array}
\right.
\renewcommand{\arraystretch}{1}
\end{equation}
The coefficients $M, M_a, \widetilde{M}_a$, and $M_{ab}$ 
are defined in a fourth-order neighborhood of a point $x \in V$ 
and satisfy the relations
 \begin{equation}\label{eq:55}
\renewcommand{\arraystretch}{1.3}
\left\{
\begin{array}{ll}
M_a - \widetilde{M}_a = K_b \lambda_a^b - 2 R^1_{11a},\\
M_{ab} = - R^1_{1ab},
\end{array}
\right.
\renewcommand{\arraystretch}{1}
\end{equation}
which are obtained if we substitute expansions (54) 
into equations (53).

For a fixed point $x \in V$, equations (54) become
 \begin{equation}\label{eq:56}
\delta K =  K \pi_1^1,
\end{equation}
and 
 \begin{equation}\label{eq:57}
\nabla_\delta K_a +  (\lambda_a^b- K \delta_a^b) \pi_b^1 = 0.
\end{equation}
Equation (56) shows that the quantity $K$ is a relative invariant 
of weight 1.

Since by theorem hypothesis, the quantity $K$ 
is not   a root of characteristic equation $(37)$, the 
affinor  
 \begin{equation}\label{eq:58}
\Lambda_a^b =   \lambda_a^b - K \delta_a^b.
\end{equation}
is nondegenerate. As the affinor $\lambda_a^b$, 
the  affinor $\Lambda_a^b$ is of weight 1. 
Thus  the inverse 
affinor $\widetilde{\Lambda}_b^a$ 
of the affinor $\lambda_b^a$ is of weight $-1$, i.e., this 
inverse affinor satisfies the equations
 \begin{equation}\label{eq:59}
\nabla_\delta 
\widetilde{\Lambda}_a^b = - \Lambda_a^b \pi_1^1.
\end{equation}

Further consider the quantities
 \begin{equation}\label{eq:60}
L_a = \widetilde{\Lambda}_a^b K_b.
\end{equation}
Differentiating equations (60) with respect to fiber parameters 
and taking into account conditions (59) and (57), we 
find that 
 \begin{equation}\label{eq:61}
\nabla_\delta L_a + L_a \pi_1^1 + \pi_a^1 = 0.
\end{equation}
Comparing equations (61) and (51), we see 
the quantities $L_a$ form a 
{\em normalizing object} of a hypersurface $V \subset (M, g)$ 
intrinsically defined by the geometry of $V$ 
and defined in its third-order neighborhood. 

Moreover, the vectors
$$
\widetilde{e}_a = e_a + L_a e_1
$$
define an invariant screen subspace $S_x$ and, 
along with it, an invariant screen distribution 
$S = \cup_{x \in V} S_x$ that is intrinsically 
connected with a lightlike hypersurface $V \subset (M, g)$. 
\rule{3mm}{3mm}

We make a  reduction in the frame bundle associated 
with a hypersurface $V$ by superposing the vectors 
$e_a$ and $\widetilde{e}_a$. Then we obtain $L_a = 0, \;
K_a = 0$, and as a result, the second group of equations 
(54) takes the form
$$
 \Lambda_a^b \omega^1_b = \widetilde{M}_a  \omega^1 + M_{ab} \omega^b.
$$
Since we assume that the tensor $\Lambda_a^b$ is 
nondegenerate, we can solve the last equations with respect 
to the 1-forms $\omega_a^1$. As a result, we obtain 
equations (29) where 
\begin{equation}\label{eq:62}
  \nu_a =  \widetilde{\Lambda}_a^b \widetilde{M}_b, \;\; 
\nu_{ab} =  \widetilde{\Lambda}_a^c M_{cb}. 
\end{equation}
These quantities are defined in a fourth-order 
neighborhood of a point $x \in V$. 
 This and (28) imply that the curvature tensor 
of the affine connection $\Gamma$ induced by 
the screen distribution $S$ we have constructed is 
defined in a fourth-order neighborhood of a point $x \in V$. 

In the same way as in Section {\bf 10}, one can prove 
that the screen distribution $S$ is integrable if and only if $\nu_{ab} = \nu_{ba}$.

Note that in the papers \cite{Be96a} and \cite{Be96b} 
as well as in the book \cite{DB96}, the authors consider 
 canonical screen distributions on a lightlike hypersurface 
$M$ of  a pseudo-Euclidean space ${\bf R}^n_q$ or 
a pseudo-Riemannian space $(\widetilde{M}, \widetilde{g})$ 
(here we used their notations).
However, this distribution and affine connections induced by them 
are not intrinsically connected 
with the geometry of a lightlike hypersurface $M$ since 
they are defined by means of a vector field $V$ connected 
with a coordinate system of the ambient space 
 ${\bf R}^n_q$ or $(\widetilde{M}, \widetilde{g})$. 
In fact, for example, 
in ${\bf R}^n_q$ this vector field $V$ is defined by 
formula (6.8) (see p. 115 of \cite{DB96}) which in the case 
$q = 1$ take the form $V = - D^0 \frac{\partial}{\partial x^0}$, i.e., 
 the vector field $V$ is a field of tangents vectors  
 to the lines $x^0$ of the curvilinear coordinate system of ${\bf R}^n_1$. 
Thus the vector field $V$ as well as the vector field $N$ (see (6.10) in [DB])
and the screen distribution $S$ (see p. 116 in [DB]) constructed by 
means of $V$ are  neither invariant nor intrinsically connected with 
the geometry of $M$. 
 
Note also that a canonical screen distribution constructed in 
\cite{Be96a}, \cite{Be96b} and \cite{DB96} is defined by elements of a 
first-order differential neighborhood of a hypersurface $M$. As we showed 
in Sections 10 and 11, screen distributions  intrinsically connected 
with the geometry of a lightlike hypersurface $M$ can be constructed 
only in a third-order differential neighborhood of $M$.

Finally note that a screen distribution similar to that in 
\cite{Be96a}, \cite{Be96b} and \cite{DB96}  
was constructed by Bonnor in 1972 (see \cite{B72}) 
who gave a physical justification for such a distribution.

{\bf 12. An affine connection on totally geodesic and totally 
umbilical lightlike hypersurfaces.} We will prove the following 
theorem. 

\begin{theorem}
The second fundamental tensor of the pseudo-Riemannian space 
$(M, g)$ vanishes on a totally geodesic lightlike hypersurface 
$V \subset (M, g)$. For any choice 
of isotropic normalization of a totally geodesic 
lightlike hypersurface $V$, an  affine connection 
is induced on $V$, and the curvature tensor 
of this connection is completely determined by 
the curvature tensor of the manifold $(M, g)$.
\end{theorem}

{\sf Proof.} 
The equations of geodesic 
lines on a pseudo-Riemannian manifold $(M, g)$ have the form 
(21). Since in a first-order frame a hypersurface $V$ 
is defined by equation (11), $V$ will be totally geodesic if 
equations (21) are identically satisfied on it.

For $i = n$, equations (21) give
$$
\omega^i \omega_i^n = 0, \;\;\;\; i = 1, \ldots , n - 1.
$$
By (17) it follows that 
$$
\lambda_{ij} = 0. 
$$

Suppose that $V \subset (M, g)$ is a lightlike hypersurface. 
In an isotropic first-order frame chosen for $V$ in Section 
{\bf 4}, the second fundamental tensor of $V$ has the form
$$
(\lambda_{ij}) = \pmatrix{0 & 0 \cr 
                         0 & \lambda_{ab} \cr}, \;\;\;\;\;
a, b = 2, \ldots, n - 1,
$$
and satisfies equation (32). If $V$ is totally geodesic, 
then on it we have
 \begin{equation}\label{eq:63}
\lambda_{ab}  = 0.
\end{equation}
From equation (34) it follows that 
$$
R^n_{ab1} = 0, \;\;\  \mu_{abc} = 0.
$$
The first of these equations shows that a totally geodesic 
lightlike hypersurface $V$ has the vanishing isotropic 
sectional curvature, $K_N (\sigma) = 0.$ Since 
the second fundamental tensor of such a $V$ also vanishes, 
it is impossible to find an invariant normalization 
of $V$ intrinsically connected with the geometry of $V$ 
by means of this tensor. 

However, an affine connection on totally geodesic 
lightlike hypersurfaces can be defined uniquely. In fact, 
equations (63) are equivalent to the equations 
$\omega_a^n = 0$. It follows from these equations that 
in structure equations (28) of the affine connection 
induced on $V$, the term $\omega_b^n \wedge \omega_n^a$ 
in the right-hand side of the last equation 
vanishes. This proves Theorem 7. \rule{3mm}{3mm}

\begin{corollary} 
If the curvature tensor of the manifold $(M, g)$ vanishes, 
$($i.e., this manifold is a Minkowski space $R^n_1$$)$, then 
totally geodesic lightlike hypersurfaces are isotropic 
hyperplanes of $R^n_1$.
\end{corollary}

Next we consider totally umbilical lightlike hypersurfaces 
$V \subset (M, g)$. They 
are defined by the equations 
 \begin{equation}\label{eq:64}
\lambda_{ab}  = \lambda g_{ab},
\end{equation}
where $\lambda \neq 0$. It follows from equations 
(64) and (25) that the isotropic geodesic $x e_1$ of 
the hypersurface $V$ carries a single singular point
 \begin{equation}\label{eq:65}
F = x - \frac{1}{\lambda} e_1.
\end{equation}
Differentiating equation (65) and applying equations (3) and 
(20), we find that 
 \begin{equation}\label{eq:66}
dF =  \frac{1}{\lambda^2} \Bigl(d\lambda - \lambda \omega_1^1 
+ \lambda^2 \omega^1\Bigr)e_1.
\end{equation}
Substituting expressions (64) into equations (34), 
we obtain that 
 \begin{equation}\label{eq:67}
g_{ab} (d\lambda - \lambda \omega_1^1 + \lambda^2 
\omega^1) + 2 R^n_{ab1} \omega^1 
= \mu_{abc} \omega^c.
\end{equation}
This implies that 
 \begin{equation}\label{eq:68}
d\lambda - \lambda \omega_1^1 + \lambda^2 
\omega^1 = \mu \omega^1 + \mu_a \omega^a. 
\end{equation}
If we substitute this expression into equations (67), 
we find that 
$$
g_{ab} (\mu \omega^1 + \mu_a \omega^a) + 2 R^n_{ab1} \omega^1 
= \mu_{abc} \omega^c.
$$
Equating coefficients in linearly independent 1-forms $\omega^1$ and $\omega^a$, we obtain
 \begin{equation}\label{eq:69}
R^n_{ab1} = - \frac{1}{2} g_{ab} \mu 
\end{equation}
and 
 \begin{equation}\label{eq:70}
g_{ab} \mu_c = \mu_{abc}. 
\end{equation}

Since the quantities $\mu_{abc}$ are symmetric with 
respect to all indices, it follows from (70) that 
$$
g_{ab} \mu_c = g_{ac} \mu_b.
$$
Contracting these equations with $g^{ab}$, we find that 
 \begin{equation}\label{eq:71}
(n - 3) \mu_c = 0. 
\end{equation}
It follows that if $n \geq 4$, then $\mu_c = 0$. Note 
that the case $n = 3$ is not interesting since 
for $n = 3$, a lightlike hypersurface 
becomes an isotropic curve. 

Now equations (68) take the form
 \begin{equation}\label{eq:72}
d\lambda - \lambda \omega_1^1 + \lambda^2 
\omega^1 = \mu \omega^1. 
\end{equation}
Taking the exterior derivative of equation (72), we find that
 \begin{equation}\label{eq:73}
(d\mu - 2 \mu \omega_1^1) \wedge \omega^1  - \mu \omega_a^1 
\wedge \omega^a + \lambda R^1_{1kl} \omega^k \wedge \omega^l = 0. 
\end{equation}
If $\lambda \neq 0$, then  for $\mu = 0$ 
equation (73) implies that 
 \begin{equation}\label{eq:74}
R^1_{1kl} = 0. 
\end{equation}
If $\mu \neq 0$, then it follows from (73) 
that 
 \begin{equation}\label{eq:75}
\renewcommand{\arraystretch}{1.3}
\left\{
\begin{array}{rl}
\displaystyle \frac{d\mu}{\mu} - 2 \omega_1^1 = \nu \omega^1 
+ \nu_a \omega^a, \\ 
- \omega^1_a = \widetilde{\nu}_a \omega^l 
+ \nu_{ab} \omega^b.
\end{array}
\right.
\renewcommand{\arraystretch}{1}
\end{equation}
Substituting these decompositions into equation (73), 
we find that 
 \begin{equation}\label{eq:76}
\nu_a - \widetilde{\nu}_a = \frac{2\lambda}{\mu} R^1_{11a}, \;\; 
\nu_{[ab]} = \frac{\lambda}{\mu} R^1_{1ab}.
\end{equation}
The quantities $\nu, \nu_a$, and $\widetilde{\nu}_a$ 
are defined in a fourth-order differential neighborhood 
of a point $x \in (M, g)$

Using equations of this section, we will 
prove further three theorems.

\begin{theorem} 
 The isotropic sectional curvature 
of a totally umbilical lightlike hypersurface $V \subset (M, g)$ 
depends on its point $x \in V$ 
and does not depend on an isotropic 
$2$-plane $\sigma = e_1 \wedge P$, where 
$P \in T_x (V)$.
\end{theorem}

{\sf Proof.} 
In fact,  it follows from (69) that 
$
R_{1ab1} = \frac{1}{2} g_{ab}\; \mu,
$
and this and formula (42) give 
$
K_N (\sigma) = \frac{1}{2}  \mu. \rule{3mm}{3mm}
$

\begin{theorem} 
 If 
for $n \geq 4$, 
the isotropic sectional curvature 
of a totally umbilical hypersurface 
 $V \subset (M, g)$ vanishes, then 
the  hypersurface $V$ is an isotropic 
 cone of the manifold $(M, g)$. 
On such a hypersurface $V$, it is impossible to 
construct an invariant normalization and an invariant 
affine connection 
intrinsically connected with the geometry of $V$. 
The components of the curvature tensor  of the manifold $(M, g)$ 
satisfy the equations
 \begin{equation}\label{eq:77}
R^n_{ab1} = 0, \;\; R^1_{1kl} = 0.
\end{equation}
\end{theorem}

{\sf Proof.}  The proof of the main 
part of this theorem follows from equations (72) 
and (66).  Since for $\mu = 0$ 
differentiation of equation (72) gives only equations 
(74) that does not contain the 
1-forms $\omega_a^1$ defining a screen distribution $S$, 
an intrinsic normalization and 
an intrinsic affine connection on such a hypersurface 
$V$ cannot be found.  
Relations (77) follows from (69) and (74). \rule{3mm}{3mm}

\begin{theorem} 
 If the isotropic sectional curvature of 
a totally umbilical manifold $(M, g)$ does not 
vanish, then a singular point $F$ of its 
isotropic geodesic $x e_1$ describes an isotropic line 
$\gamma$. On $V$ one can  define an invariant 
 screen distribution $S$ intrinsically connected with 
the geometry of $V$. This distribution is 
integrable if and only if $R^1_{1ab} = 0$.
\end{theorem}

{\sf Proof.} In fact, by (66) and (72), we have 
 \begin{equation}\label{eq:78}
dF = \frac{\mu}{\lambda^2} \omega^1 e_1.
\end{equation}
This means that the point $F$ describes a line $\gamma$ 
tangent to the vector $e_1$, i. e., an isotropic curve. 
The equation $\omega^1 = 0$ defines on $V$ a 
 screen distribution $S$ intrinsically connected with 
the geometry of $V$. If a point $x$ moves 
along integral lines of the distribution $S$, then by (78) 
the point $F$ is fixed. It follows from the second 
equation of (76) that the screen distribution $S$ 
is integrable if and only if the components $R^1_{1ab}$ 
of the curvature tensor of the manifold 
$(M, g)$ vanish on $V$, $R^1_{1ab} = 0$. In this case 
the fibration of isotropic geodesics is decomposed 
into a one-parameter family of cones.
\rule{3mm}{3mm}

{\bf 13. Lightlike hypersurfaces on a pseudo-Riemannian manifold 
$(M, g)$ of Lorentzian signature and constant curvature.}

The tensor of Riemannian curvature of a Riemannian 
or pseudo-Riemannian manifold $(M, g)$ of constant curvature 
has the form 
 \begin{equation}\label{eq:79}
R_{ijkl} = K (g_{ik} g_{jl} - g_{il} g_{jk}),
\end{equation}
where $K$ is the curvature of the manifold. By Schur's theorem 
(see \cite{S86} or \cite{KN63}, pp. Section {\bf 5.3}), 
for $n \geq 3$, 
the curvature $K$ does not depend on a point $x \in (M, g)$, i.e, 
$K$ is constant on the manifold $(M, g)$. 

For $K = 0$,  the manifold 
$(M, g)$ of   Lorentzian signature and constant curvature is 
the Minkowski space $R^n_1$; for $K > 0$, it is 
the de Sitter space $S^n_1$ of first kind whose 
projective model was considered in detail in \cite{AG98a} 
and \cite{AG98c}; and for $K < 0$, it is 
the de Sitter space $H^n_1$ of second kind (see \cite{BEE96}, 
pp. 115--117).

Harris in \cite{H82} proved the following theorem. 

\begin{theorem} 
A pseudo-Riemannian manifold 
$(M, g)$ of  Lorentzian signature has 
a  constant curvature if and only if  
its  isotropic sectional curvature $K_N (\sigma)$ 
vanishes.
\end{theorem}

{\sf Proof.} 
It is not so difficult to prove 
the necessity  of this theorem. In fact, consider 
an isotropic frame bundle   on 
a manifold $(M, g)$.  In this frame bundle the metric tensor 
$g_{ij}$ has the form (8). This and equations (79) imply 
that 
 \begin{equation}\label{eq:80}
R_{1ab1} = 0.
\end{equation}
But since $e_1$ is an arbitrary isotropic vector, 
by (43), condition (80) means that  $K_N (\sigma) = 0$ 
on the manifold $(M, g)$. 

The proof of sufficiency is more complicated (see \cite{H82}). 
\rule{3mm}{3mm}

By conditions (80) equations (34) on a lightlike 
hypersurface on a manifold $(M, g)$ of  constant curvature 
take the form
\begin{equation}\label{eq:81}
\nabla \lambda_{ab} - \lambda_{ab} \omega_1^1 
+ \lambda _{ac} g^{ce} \lambda_{eb} \omega^1 
 = \mu_{abc} \omega^c. 
\end{equation}
As a result the covariant derivative of the 
tensor $\lambda_{ab}$ in the direction of the 
vector $e_1$ has the following expression:
$$
(\nabla \lambda_{ab} - \lambda_{ab} \omega_1^1)_{,1} 
= - \lambda _{ac} g^{ce} \lambda_{eb} 
$$
It is expressed only in terms of quantities 
defined in a second-order differential neighborhood 
of a point $x \in (M, g)$.

A construction of an invariant normalization 
and an invariant affine connection for a 
lightlike hypersurface $V \subset (M, g)$ of constant 
curvature can be done in the same way as 
in general case following the scheme indicated 
in Sections {\bf 10} and {\bf 11} with the only difference 
that in formulas (46) and (52) 
the quantity $K$ is  
defined now in a second-order differential neighborhood 
of a point $x \in (M, g)$ (not the third-order as this was 
in the general case). 

Consider a totally umbilical lightlike hypersurface $V$ 
on a manifold $(M, g)$ of Lorentzian signature and constant 
curvature. 
 By Theorem 12, on such a hypersurface 
the  isotropic sectional curvature $K_N (\sigma)$ 
vanishes. This and Theorem 10 imply the following result.

\begin{theorem} 
Totally umbilical lightlike hypersurface $V$ 
on a manifold $(M, g)$ of  Lorentzian 
signature and constant curvature are the light cones 
of $(M, g)$.
\end{theorem}

Note that any Riemannian or pseudo-Euclidean manifold 
$(M, g)$ of constant curvature is conformally flat (see 
for example, \cite{R59}, \S 122). 
Hence Theorem 13 follows from Theorem 7, part b, 
 of the paper \cite{AG98c}.

{\bf 14. An intrinsic normalization of a lightlike hypersurface 
$V$ on a four-dimensional manifold $(M, g)$ of Lorentzian 
signature.} 
Consider a lightlike hypersurface on a manifold $(M, g), \; 
\dim\; M = 4, \;
{{\rm sign}}\; g = (3, 1)$. All formulas of Sections {\bf 4-8} 
hold on such a hypersurface, and the range of the indices $a, b, c$ 
is $2, 3$: $a, b, c = 2, 3$. 
We reduce simultaneously the first and the second fundamental 
tensors of the hypersurface $V$ to diagonal forms:
\begin{equation}\label{eq:82}
(g_{ab}) = \pmatrix{1 & 0 \cr
                  0 & 1 \cr}, \;\;
(\lambda_{ab}) = \pmatrix{\lambda_2 & 0 \cr 
                         0 & \lambda_3}
\end{equation}
and  assume that $\displaystyle \frac{\lambda_2}{\lambda_3} 
\neq {{\rm const}}$, and $\lambda_2 \neq 0, \lambda_3 \neq 0$.

From the last equation of (10) and the first relation of (82) it 
follows that on $V$ we have 
\begin{equation}\label{eq:83}
\omega_2^2 = \omega_3^3 = 0, \;\;\omega_2^3 + \omega_3^2 = 0,
 \end{equation}
and equations (34) take the form 
\begin{equation}\label{eq:84}
\renewcommand{\arraystretch}{1.3}
\left\{
\begin{array}{ll}
d\lambda_2 - \lambda_2 \omega_1^1 + ((\lambda_2)^2 + R^4_{221})\omega^1 
= \mu_{22c} \omega^c, \\
d\lambda_3 - \lambda_3 \omega_1^1 + ((\lambda_3)^2 + R^4_{331})\omega^1 = \mu_{32c} \omega^c, \\
(\lambda_2 - \lambda_3) \omega_2^3 + R^4_{231}\omega^1 
= \mu_{23c} \omega^c.
\end{array}
\right.
\renewcommand{\arraystretch}{1}
\end{equation}
Since $\lambda_2 \neq \lambda_3$, then the last equation implies 
that 
\begin{equation}\label{eq:85}
\omega_2^3 = \displaystyle \frac{1}{\lambda_2 - \lambda_3} 
(R_{1231} \omega^1 + \mu_{232} \omega^2 + \mu_{233} \omega^3).
 \end{equation}
The first two equations of (84) can be written as
\begin{equation}\label{eq:86}
\renewcommand{\arraystretch}{1.3}
\left\{
\begin{array}{ll}
d \lambda_2 
- \lambda_2 \omega_1^1 
=(R_{121}- (\lambda_2)^2) \omega^1 + \mu_{222} \omega^2 + \mu_{223} \omega_3, \\
d\lambda_3 - \lambda_3 \omega_1^1 = (R_{1331} - (\lambda_3)^2)\omega^1 
+ \mu_{332} \omega^2 + \mu_{333} \omega^3.
\end{array}
\right.
\renewcommand{\arraystretch}{1}
\end{equation}

The quantities $\lambda_2$ and $\lambda_3$ are relative invariants of weight one. The equations to which these invariants satisfy can be written in the form (52), where
\begin{equation}\label{eq:87}
\renewcommand{\arraystretch}{1.3}
\left\{
\begin{array}{lll}
K_2 = \lambda_2 - \displaystyle \frac{R_{1221}}{\lambda_2}, &
K_{22} =  - \displaystyle \frac{\mu_{222}}{\lambda_{23}}, &
K_{23} = - \displaystyle \frac{\mu_{223}}{\lambda_2}, \\
K_3 = \lambda_3 - \displaystyle \frac{R_{1331}}{\lambda_3}, &
K_{32} = - \displaystyle \frac{\mu_{332}}{\lambda_3}, &
K_{33} = - \displaystyle \frac{\mu_{333}}{\lambda_3}. 
\end{array}
\right.
\renewcommand{\arraystretch}{1}
\end{equation}
The first index in these equations is the index of the 
relative invariant $\lambda_a$. 

By Theorem 6, if the coefficients $K_a$ are not roots of the 
characteristic 
equation of the affinor $(\lambda_b^a)$, then by means of the 
coefficients 
$K_{ab}$ we can construct the normalizing objects $L_{ab}$. These normalizing objects determine two invariant normalizations 
intrinsically connected with the 
geometry of the hypersurface $V$.

The ratio $\displaystyle \frac{\lambda_2}{\lambda_3}$ of the 
eigenvalues of 
the affinor $(\lambda_b^a)$ is an absolute invariant. It follows 
from 
equations (86) that this absolute invariant satisfies the 
equation
\begin{equation}\label{eq:88} 
\ln \displaystyle \Bigl|\frac{\lambda_2}{\lambda_3}\Bigr| 
= \Biggl(\displaystyle \frac{K_2}{\lambda_2}- \frac{K_3}{\lambda_3}\Biggr)\omega^1  
+ \Biggl(\displaystyle \frac{K_{22}}{\lambda_2}- \frac{K_{32}}{\lambda_3}\Biggr)\omega^2 + \Biggl(\displaystyle \frac{K_{23}}{\lambda_2}- 
\frac{K_{33}}{\lambda_3}\Biggr)\omega^3.
 \end{equation}
By Theorem 5, if the coefficient in $\omega^1$ in equation (88)
 is different from 0 (i.e., if the quantities $K_2$ and $K_3$ are 
 not proportional to the eigenvalues $\lambda_2$ and $\lambda_3$ 
of  the affinor $(\lambda_b^a)$), then the absolute invariant 
$\displaystyle \frac{\lambda_2}{\lambda_3}$ allows us to 
construct 
one more invariant normalization intrinsically connected with 
the geometry of the hypersurface $V$. The screen distribution 
defining this normalization is tangent to level submanifolds 
of the invariant $\displaystyle \frac{\lambda_2}{\lambda_3}$.

Thus we have proved the following result.

\begin{theorem} If the eigenvalues $\lambda_2$ and $\lambda_3$ of 
the affinor $(\lambda_b^a)$ of a lightlike hypersurface 
$V \subset (M, g), \; \dim M = 4$, are different from $0$, 
  the absolute invariant 
$\displaystyle \frac{\lambda_2}{\lambda_3} \neq {{\rm const}}$, 
and the coefficients $K_2$ and $K_3$ defined by formulas $(87)$ do not coincide with any of the eigenvalues   $\lambda_2$ and $\lambda_3$ and are not proportional to them, then on such a 
hypersurface we can construct three 
invariant normalizations intrinsically connected with the
 geometry of $V$, and the screen distribution of one of 
these normalizations is integrable.
\end{theorem}

Note also that the eigenvectors $e_2$ and $e_3$ 
corresponding to the eigenvalues $\lambda_2$ and $\lambda_3$ of 
the affinor $(\lambda_b^a)$ generate  two orthogonal vector 
fields 
on screen distributions of normalizations we have constructed. 
These vector fields with the field of isotropic vectors $e_1$ 
determine the coordinate net on the hypersurface $V$. In general, 
this net is 
not holonomic. This  the means that in general,  the    
two-dimensional distributions defined by the 
eigenvectors of the affinor $(\lambda_b^a)$ and the vectors $e_1$ 
are not integrable.

\noindent
{\em Authors' addresses}:\\

\noindent
\begin{tabular}{ll}
M.A. Akivis &                           V.V. Goldberg \\
Department of Mathematics      &      Department of Mathematics\\
Jerusalem College of Technology  & 
                       New Jersey Institute of Technology \\
- Mahon Lev, P. O. B. 16031 & University Heights \\
Jerusalem 91160, Israel & Newark, NJ 07102, U.S.A.  \\
&\\
E-mail address: akivis@avoda.jct.ac.il &E-mail address: 
vlgold@numerics.njit.edu
\end{tabular}
\end{document}